\documentclass{amsart}
\usepackage{graphicx}
\newtheorem{thm}{Theorem}[section]

\newtheorem{lem}[thm]{Lemma}
\newtheorem{prop}[thm]{Proposition}
\theoremstyle{definition}

\theoremstyle{remark}

\numberwithin{equation}{section}
\begin{document}

\title[Special exponential sums]{A note on some special exponential sums}%
\author{Xiwang Cao}%
\address{Xiwang Cao is with the School of Mathematical Sciences, Nanjing University of
Aeronautics and Astronautics, Nanjing 210016, China, email: {\tt
xwcao@nuaa.edu.cn},
%Lei Hu is with the State Key Laboratory of Information security, Institute of Information Engineering, Chinese Academy of Sciences, Beijing 100093, China.
%email:{\tt hu@is.ac.cn
}%
%\thanks{}%

%\date{}%
%\dedicatory{}%
%\commby{}%
% ----------------------------------------------------------------
\begin{abstract}
In this note, we presented a new decomposition of elements of finite fields of even order and illustrated that it is an effective tool in evaluation of some specific exponential sums over finite fields, the explicit value of some exponential sums were obtained.
\end{abstract}
\maketitle
% ----------------------------------------------------------------

\section{Introduction}

Let $n=2m$ be a positive even integer and $\mathbb{F}_{2^n}$ the finite field of order $2^n$. Let ${\rm Tr}_1^n(x)=\sum_{i=0}^{n-1}x^{2^i}$ be the trace function. We encountered the following exponential sums when we investigate the Walsh spectra of some boolean functions.
\begin{eqnarray}\label{f-1}
   && p(\mu)=\sum_{a\in \mathbb{F}_{2^n}\setminus \mathbb{F}_2}\chi_n\left(\mu\frac{a^{2^m}+a}{a^2+a}\right), {\rm and}\\ \label{f-1'}&&q(\mu)=\sum_{a\in \mathbb{F}_{2^n}\setminus \mathbb{F}_{2^m}}\chi_n\left(\mu\frac{a^2+a}{a^{2^m}+a}\right),\\
   \label{f-1''}&&q_{s}(\mu)=\sum_{a\in \mathbb{F}_{2^n}\setminus \mathbb{F}_{2^m}}\chi_n\left(\mu\frac{(a^2+a)^{2^s}}{a^{2^m}+a}\right),\\
   \label{f-l}&&r(l)=\sum_{a\in \mathbb{F}_{2^n}}\chi_n((a^{2^m}+a)L(a)),
\end{eqnarray}
where $\mu\in \mathbb{F}_{2^m}$, $\chi_n(x)=(-1)^{{\rm Tr}_1^n(x)}$ and $s$ is a positive integer with $\gcd(s,m)=d$. $L(x)=\sum_{i=0}^k\alpha_ix^{2^{a_i}}\in \mathbb{F}_{2^m}[x]$ is a linearized polynomial with coefficients in $\mathbb{F}_{2^m}$.

These
sums form a subset of a much larger class of exponential sums of the form
\begin{equation}\label{f-2}
    \sum_{x\in \mathbb{F}_q}\chi_n(f(x))
\end{equation}
where $f(x)\in \mathbb{F}_q[X]$. The sums of the form (\ref{f-2}) are also known as Weil sums. The problem of
explicitly evaluating Weil sums is quite often difficult. Results giving estimates
for the absolute value of the Weil sum are more common and such results have been
regularly appearing for many years. We refer the reader to \cite{lidl} for an
overview of the related researches.

%Let $r,s$ be two integers. We define two more exponential sums as follows.
%\begin{equation}\label{f-3}
%   p_{r,s}(\mu)=\sum_{a\in \mathbb{F}_{2^n}\setminus \mathbb{F}_2}\chi_n\left(\mu\frac{(a^{2^m}+a)^r}{(a^2+a)^{2^s}}\right), {\rm and}\  q_{r,s}(\mu)=\sum_{a\in \mathbb{F}_{2^n}\setminus \mathbb{F}_{2^m}}\chi_n\left(\mu\frac{(a^2+a)^{2^s}}{(a^{2^m}+a)^r}\right).
%\end{equation}
%When $s=0$ and $r=1$, the sums in (\ref{f-3}) reduce to that in (\ref{f-1}).
In this note, we will explicitly evaluate these exponential sums. The main results are the following:
{\thm \label{thm01}For every $\mu\in \mathbb{F}_{2^m}^*$, one has that
\begin{eqnarray}\label{f-4}
   && p(\mu)=\sum_{a\in \mathbb{F}_{2^n}\setminus \mathbb{F}_2}\chi_n\left(\mu\frac{a^{2^m}+a}{a^2+a}\right)=-2-(1+k_m(\mu))^2,\\
   &&q(\mu)=\sum_{a\in \mathbb{F}_{2^n}\setminus \mathbb{F}_{2^m}}\chi_n\left(\mu\frac{a^2+a}{a^{2^m}+a}\right)=-2^m\chi_m(\mu).
\end{eqnarray}}
{\thm If $\gcd(s,m)=d=1$, and $m$ is odd, then

(1) $q_s(\mu)=-2^m$ if and only if ${\rm Tr}_1^m(\mu^{\frac{1}{2^s+1}})=0$;

(2) if ${\rm Tr}_1^m(\mu^{\frac{1}{2^s+1}})=1$, then there is an $h\in L$ such that  $\mu^{\frac{1}{2^s+1}}=h^{2^{s}}+h^{2^{m-s}}+1$ and

\begin{equation*}
    q_s(\mu)=2^m\left(\chi_m(h^{2^{s}+1}+h)\left(\frac{2}{m}\right)2^{(m+1)/2}-1\right),
\end{equation*}
where $\left(\frac{2}{m}\right)$ is the Jacobi symbol, and $k_m(\mu)$ is the Kloosterman sum.
}
{\thm For every $L(x)\sum_{i=0}^k\alpha_ix^{2^{a_i}}\in \mathbb{F}_{2^m}[x]$, we have
\begin{equation}\label{f-1102}
    r(l)=2^m\sum_{u\in \mathbb{F}_{2^m}}\chi_m(uL(u))=2^m\sum_{u\in \mathbb{F}_{2^m}}\chi_m\left(\sum_{i=0}^k\alpha_iu^{2^{a_i}+1}\right).
\end{equation}}

\section{Notations and Preliminaries}

\subsection{ Trace representations of Boolean functions}

Let $n$ be a positive integer and $\mathbb{F}_{2^n}$ be the finite field with $2^n$ elements. A Boolean function on $\mathbb{F}_{2^n}$ is
an $\{0,1\}$-valued function from $\mathbb{F}_{2^n}$ to $\mathbb{F}_2$.

For any positive integer $n$, and for any positive integer $k$ dividing $n$, the trace
function from $\mathbb{F}_{2^n}$ to $\mathbb{F}_{2^k}$, denoted by ${\rm Tr}_k^n$, is the mapping defined as
\begin{equation*}
   {\rm Tr}_k^n(x)=x+x^{2^k}+x^{2^{2k}}+\cdots+x^{2^{n-k}}.
\end{equation*}
In particular, the absolute trace over $\mathbb{F}_2$ is the function
${\rm Tr}_1^n(x)=\sum_{i=0}^{n-1}x^{2^i}$ for $k=1$. Recall that, for every integer $k$ dividing
$n$, the trace function satisfies the transitivity property, that
is, for all $x\in \mathbb{F}_{2^n}$, it holds that \cite{lidl}
\begin{equation*}
    {\rm Tr}_1^n(x)={\rm Tr}_1^k({\rm Tr}_k^n(x)).
\end{equation*}
It is known \cite{carlet} that every nonzero Boolean function $g$ defined on $\mathbb{F}_{2^n}$ has a
unique trace expansion of the form
\begin{equation*}
    g(x)=\sum_{j\in \Gamma_n}{\rm Tr}_1^{o(j)}(ax^j)+\epsilon(1+x^{2^n-1})
\end{equation*}
where $\Gamma_n$ is the set of integers obtained by choosing one element in each cyclotomic coset of $2$ modulo $(2^n-1)$ and $o(j)$ is the size of the cyclotomic coset of $2$ modulo $(2^n-1)$  containing $j$, $a_j\in \mathbb{F}_{2^{o(j)}}$, and $\epsilon=wt(g)$ modulo $2$ where $wt(g)$ is the Hamming weight
of the image vector of $g$, that is, the cardinality of its support set
\begin{equation*}
Supp(g):=\{x\in \mathbb{F}_{2^n}|g(x)=1\}.
\end{equation*}

\subsection{Polar decomposition}

Let $n=2m$. Denote the subgroup of $(2^m+1)$-th roots of unity in $\mathbb{F}_{2^n}$ by $\mathfrak{S}$, i.e., $\mathfrak{S}=\{z\in \mathbb{F}_{2^n}| z^{2^m+1}=1\}$. For every $x\in \mathbb{F}_{2^n}^*=\mathbb{F}_q\setminus \{0\}$, there is a unique polar decomposition of $x$ as $x=yz$ where $y\in \mathbb{F}_{2^m}^*$ and $z\in \mathfrak{S}$. In fact, $y=x^{(2^m+1)2^{m-1}}, z=x^{(2^m-1)2^{m-1}}$. Denote $x^{2^m}$ by $\overline{x}$. Then for every $x\in \mathbb{F}_{2^n}^*$, $x\in \mathbb{F}_{2^m}$ if and only if $x=\overline{x}$, and $x\in \mathfrak{S}$ if and only if $\overline{x}=x^{-1}$. it is evident that for every $x\in \mathbb{F}_{2^n}^*$, one has that $x+\overline{x}, x\overline{x}\in \mathbb{F}_{2^m}$ and $x/\overline{x},\overline{x}/x\in \mathfrak{S}$. Note that $x\mapsto \overline{x}$ is an isomorphism of the finite field $\mathbb{F}_{2^n}$.

\subsection{ Kloosterman sums}

For every $a,b\in \mathbb{F}_{2^m}$, Kloosterman sum is defined by
\begin{eqnarray}
  \label{equ-k}k_m(a,b) &=& \sum_{x\in \mathbf{F}_{2^m}^*}(-1)^{{\rm Tr}_m(ax+bx^{-1})}.
 \end{eqnarray}
 It is easy to check that $k_m(a,b)=k_m(ab,1)=k_m(1,ab)$. For simplicity, denote $k_m(a,1)=k_m(1,a)$ by $k_m(a)$. Moreover, The Kloosterman sum $k_m(a,b)$ can be calculated recursively, that is,
 if we define
 \begin{equation*}
    k_m^{(s)}(a)=\sum_{\gamma\in \mathbb{F}_{{2^{ms}}}^*}\chi^{(s)}(a\gamma+\gamma^{-1}), a\in \mathbb{F}_{2^m},
 \end{equation*}
where $\chi^{(s)}$ is the lifting of $\chi(x)=(-1)^{{\rm Tr}(x)}$ to $\mathbb{F}_{2^{ms}}$, then
 \begin{equation}\label{f-41}
    k_m^{(s)}(a)=-k_{m}^{(s-1)}(a)k_m^{(1)}(a)-2^mk_{m}^{(s-2)}(a),
 \end{equation}
 where we put $k_m^{(0)}(a,b)=-2$ and $k_m^{(1)}(a)=k(a)$. Moreover, for all $a,b\in \mathbb{F}_{2^m}^*$, one has that
  \begin{equation*}
   |k_m(a,b)|\leq 2\sqrt{2^m},
\end{equation*}
See \cite{lidl} for details.
% Using the recursive relation, one can obtain the formulae for $k_m(1)$ as follows:
 %{\lem \cite{lidl} $k_m(1)=-\omega_1^m-\omega_2^m$, where $\omega_1,\omega_2$ are the complex roots of the equation $x^2+x+2=0$. Moreover, one has that $k_1(1)=1,k_2(1)=3$ and
 % \begin{equation}\label{f-43}
%    k_{m+2}(1)+k_{m+1}(1)+2k_{m}(1)=0.
%\end{equation}}

% {\lem \cite{carlitz} The value of $k_m(1)$ is
%\begin{eqnarray*}
 % k_m(1)&=& -\sum_{j=0}^{\lfloor m/2\rfloor}(-1)^{m-j}\frac{m}{m-j}{{m-j}\choose{j}}2^j.
%\end{eqnarray*}}

Note also that the values of Kloosterman sums over $\mathbb{F}_{2^m}$
were determined by Lachaud and Wolfmann in \cite{lachaud}.

\begin{lem}\label{lem-43} (\cite{lachaud}) The set $\{k_m(\lambda),\lambda\in \mathbb{F}_{2^m}\}$ is the set of all the integers $s\equiv -1({\rm mod}\ 4)$ in the range
\begin{equation*}
    \left[-2^{\frac{m}{2}+1},2^{\frac{m}{2}+1}\right].
\end{equation*}
\end{lem}
\subsection{A new decomposition of elements of $\mathbb{F}_{2^n}$ related to an affine subspace}
Let $n=2m$. Denote $$E=\{\lambda\in \mathbb{F}_{2^n}: \lambda^{2^m}+\lambda=1\}.$$
Then $E$ is an affine subspace of $\mathbb{F}_{2^n}/\mathbb{F}_2$. For every $x\in \mathbb{F}_{2^n}^*\setminus \mathbb{F}_{2^m}$, there is a unique pair $(u,\lambda)\in \mathbb{F}^*_{2^m}\times E$ such that $x=u\lambda$. If $x\in \mathbb{F}_{2^m}$, we just write $x=u$. This decomposition is unique, for if there are $u_1,u_2\in \mathbb{F}_{2^m}^*$ and $\lambda_1,\lambda_2\in E$ satisfying $u_1\lambda_1=u_2\lambda_2$, then
\begin{equation*}
    1=\overline{\lambda_1}+\lambda_1=(u_2/u_1)(\overline{\lambda_2}+\lambda_2)=u_2/u_1
\end{equation*}
which implies that $u_1=u_2$ and $\lambda_1=\lambda_2$.
Under this decomposition, the following two facts are easily verified.

{\bf Fact }(i) For every $x\in \mathbb{F}_{2^n}\setminus \mathbb{F}_{2^m}$, let $x=u\lambda$, $u\in \mathbb{F}_{2^m}^*, \lambda\in E$. Then ${\rm Tr}_1^n(x)={\rm Tr}_1^m(u)$.

{\bf  Fact }(ii) The map $\sigma: E \rightarrow \mathbb{F}_{2^m}; \lambda\mapsto \lambda\overline{\lambda}$ is a two-to-one map, the image set is precisely the set of elements in $\mathbb{F}_{2^m}$ which is of trace one.

\section{proof of the main results}
\subsection{The evaluation of $p(\mu)$}
For the first exponential sums, we compute that
\begin{eqnarray*}
  p(\mu)&=&\sum_{a\in \mathbb{F}_{2^n}\setminus \mathbb{F}_2}\chi_n\left(\mu\frac{a^{2^m}+a}{a^2+a}\right) \\
   &=&2^m-2+ \sum_{a\in \mathbb{F}_{2^n}\setminus \mathbb{F}_{2^m}}\chi_n\left(\mu\frac{a^{2^m}+a}{a^2+a}\right).
\end{eqnarray*}
For every $a\in \mathbb{F}_{2^n}\setminus \mathbb{F}_{2^m}$, let $a=u\lambda$ where $u\in \mathbb{F}_{2^m}^*,\lambda\in E$. Then
\begin{eqnarray*}
  &&\sum_{a\in \mathbb{F}_{2^n}\setminus \mathbb{F}_{2^m}}\chi_n\left(\mu\frac{a^{2^m}+a}{a^2+a}\right) \\
  &=&\sum_{u\in \mathbb{F}_{2^m}^*,\lambda\in E}\chi_n\left(\mu (\frac{1}{\lambda}+\frac{1}{\lambda+u})\right)\\
%&=&\sum_{u\in \mathbb{F}_{2^m}^*,\lambda\in E}\chi_m\left(\mu u^{r-2^{s}}(\frac{1}{\overline{\lambda}\lambda}+\frac{1}{(\lambda+u^{-2^s})(\overline{\lambda}+u^{-2^s})})\right)\\
&=&\sum_{u\in \mathbb{F}_{2^m}^*,\lambda\in E}\chi_m\left(\mu (\frac{1}{\overline{\lambda}\lambda}+\frac{1}{\lambda\overline{\lambda}+u^2+u})\right).
\end{eqnarray*}
By the Fact (ii), one has that
\begin{eqnarray*}
  &&\sum_{u\in \mathbb{F}_{2^m}^*,\lambda\in E}\chi_m\left(\mu (\frac{1}{\overline{\lambda}\lambda}+\frac{1}{\lambda\overline{\lambda}+u^{2}+u})\right) \\
  &=&2\sum_{v\in \mathbb{F}_{2^m}, {\rm Tr}_1^m(v)=1}\chi_m(\frac{\mu }{v})\sum_{u\in \mathbb{F}_{2^m}^*}\chi_m\left(\frac{\mu}{v+u^2+u}\right)\\
  &=&2\sum_{v\in \mathbb{F}_{2^m}, {\rm Tr}_1^m(v)=1}\chi_m(\frac{\mu }{v})\left(\chi_m(\mu/v)+2\sum_{u\in \mathbb{F}_{2^m}^*\setminus\{v\}, {\rm Tr}_1^m(u)=1}\chi_m(\frac{\mu}{u})\right).
\end{eqnarray*}
Since
\begin{eqnarray*}
 &&\chi_m(\mu/v)+2\sum_{u\in \mathbb{F}_{2^m}^*\setminus\{v\}, {\rm Tr}_1^m(u)=1}\chi_m(\frac{\mu}{u}) \\
 &=&\chi_m(\mu/v)-2\chi_m(\mu/v)+2\sum_{u\in \mathbb{F}_{2^m}, {\rm Tr}_1^m(u)=1}\chi_m(\mu/u)\\
 &=&\chi_m(\mu/v)-2\chi_m(\mu/v)+\sum_{u\in \mathbb{F}_{2^m}}\chi_m(\mu/u)(1-\chi_m(u))\\
 &=&-\chi_m(\mu/v)-\sum_{u\in \mathbb{F}_{2^m}}\chi_m(u+\mu/u)\\
 &=&-1-\chi_m(\mu/v)-k_m(\mu),
\end{eqnarray*}
one has that
\begin{eqnarray*}
 &&\sum_{a\in \mathbb{F}_{2^n}\setminus \mathbb{F}_{2^m}}\chi_n\left(\mu\frac{a^{2^m}+a}{a^2+a}\right)  \\
  &=& -2(1+k_m(\mu))\sum_{v\in \mathbb{F}_{2^m}, {\rm Tr}_1^m(v)=1}\chi_m(\mu/v)-2\sum_{v\in \mathbb{F}_{2^m}, {\rm Tr}_1^m(v)=1}1\\
  &=&-(1+k_m(\mu))^2-2^m.
\end{eqnarray*}
Therefore, we have that
\begin{equation}\label{f-4}
    p(\mu)=\sum_{a\in \mathbb{F}_{2^n}\setminus \mathbb{F}_2}\chi_n\left(\mu\frac{a^{2^m}+a}{a^2+a}\right)=-2-(1+k_m(\mu))^2.
\end{equation}
\subsection{The evaluation of $q(\mu)$}
It is evident that
\begin{eqnarray*}
 &&q(\mu)=\sum_{a\in \mathbb{F}_{2^n}\setminus \mathbb{F}_{2^m}}\chi_n(\mu\frac{a^2+a}{a^{2^m}+a}) \\
  &=& \sum_{u\in \mathbb{F}_{2^m}^*,\lambda\in E}\chi_n(\mu(u\lambda^2+\lambda))\\
  &=&\sum_{u\in \mathbb{F}_{2^m}^*,\lambda\in E}\chi_m(\mu(u+1))\\
  &=&-2^m\chi_m(\mu).
\end{eqnarray*}
\subsection{The evaluation of $q_s(\mu)$}
It is easy to see that
\begin{eqnarray*}
 &&q_s(\mu)=\sum_{a\in \mathbb{F}_{2^n}\setminus \mathbb{F}_{2^m}}\chi_n(\mu\frac{(a^2+a)^{2^s}}{a^{2^m}+a}) \\
  &=& \sum_{u\in \mathbb{F}_{2^m}^*,\lambda\in E}\chi_n(\mu (u^{2^{s+1}-1}\lambda^{2^{s+1}}+u^{2^s-1}\lambda^{2^s}))\\
  &=&\sum_{u\in \mathbb{F}_{2^m}^*,\lambda\in E}\chi_m(\mu(u^{2^{s+1}-1}+u^{2^s-1}))\\
  &=&2^m\sum_{u\in \mathbb{F}_{2^m}^*}\chi_m(\mu(u^{(2^{s}-1)(2^s+1)}+u^{2^s-1})).
\end{eqnarray*}

%{\bf Case (1)}.
Suppose that $m$ is odd and $\gcd(s,m)=d=1$.
In this case, it is obvious that $\gcd(2^s+1,2^m-1)=\gcd(2^s-1,2^m-1)=1$ and we know that \begin{eqnarray*}
 &&q_s(\mu)=\sum_{a\in \mathbb{F}_{2^n}\setminus \mathbb{F}_{2^m}}\chi_n(\mu\frac{(a^2+a)^{2^s}}{a^{2^m}+a}) \\
    &=&2^m\sum_{u\in \mathbb{F}_{2^m}^*}\chi_m(\mu(u^{(2^{s}-1)(2^s+1)}+u^{2^s-1}))\\
    &=&2^m\sum_{u\in \mathbb{F}_{2^m}^*}\chi_m(\mu(u^{2^s+1}+u)).
\end{eqnarray*}
The exponential sum $\sum_{u\in \mathbb{F}_{2^m}^*}\chi_m(\mu(u^{2^s+1}+u))$ is a special case of the following exponential sum:
\begin{equation}
   \label{equ-c} C_m^{(s)}(a,b) =\sum_{x\in \mathbf{F}_{2^m}}\chi_m(ax^{2^s+1}+bx),a,b\in \mathbb{F}_{2^m}.
\end{equation}
For odd $m$ and $\gcd(s,m)=1$, it is proved in \cite{lmw} that
{\lem \label{LMW} If $m$ is odd and
$\gcd(s,m)=1$, then

\begin{equation}\label{f-11022}
    C_m^{(s)}(1,1)=\left(\frac{2}{m}\right)2^{(m+1)/2}=\left\{\begin{array}{cc}
      2^{(m+1)/2}, & \mbox{ if } m\equiv \pm 1 ({\rm mod} \ 8), \\
      -2^{(m+1)/2}, & \mbox{ if } m\equiv \pm 3 ({\rm mod} \ 8), \\
    \end{array}\right.
\end{equation}
where $\left(\frac{2}{m}\right)$ is the Jacobi symbol.}

If $m$ is odd and $\gcd(s,m)=1$, then $x\mapsto x^{2^s+1}$ is a permutation on $L=:\mathbb{F}_{2^m}$. In this case, we have

{\prop \label{lem-sb}If $m$ is odd and $\gcd(s,m)=1$, then

(1) $C_m^{(s)}(a,b)=C_m^{(s)}(1,b/a^{\frac{1}{2^s+1}})$;

(2) $C_m^{(s)}(1,a)=C_m^{(s)}(1,a^2)$ for all $a\in L$;

(3) $C_m^{(s)}(1,a)=0$ if and only if ${\rm Tr}(a)=0$;

(4) if ${\rm Tr}(a)=1$, then there is an $h\in L$ such that  $a=h^{2^{s}}+h^{2^{m-s}}+1$ and
$$C_m^{(s)}(1,a)=\chi_m(h^{2^{s}+1}+h)C_m^{(s)}(1,1)=\chi_m(h^{2^{s}+1}+h)\left(\frac{2}{m}\right)2^{(m+1)/2}.$$

(5) Let $M_{+}=|\{a\in L|C_m^{(s)}(1,a)=\left(\frac{2}{m}\right)2^{(m+1)/2},\}|, M_{-}=|\{a\in L|C_m^{(s)}(1,a)=-\left(\frac{2}{m}\right)2^{(m+1)/2}\}|$. Then
\begin{equation}\label{f-9270}
    2M_{+}=|\{h\in L|{\rm Tr}(h^{2^sk+1}+h)=0\}|,  2M_{-}=|\{h\in L|{\rm Tr}(h^{2^s+1}+h)=1\}|,
\end{equation}
and
\begin{equation}\label{f-927}
    M_{+}=2^{m-2}+(-1)^{\frac{m^2-1}{8}}2^{\frac{m-3}{2}}, M_{-}=2^{m-2}-(-1)^{\frac{m^2-1}{8}}2^{\frac{m-3}{2}}.
\end{equation}

(6) Let
\begin{equation*}
    \mathcal{N}_{i,j}=\{h\in \mathbb{F}_{2^m}|{\rm Tr}(h^{2^s+1})=i,{\rm Tr}(h)=j\}, i,j=0,1.
\end{equation*}
And $N_{i,j}=|\mathcal{N}_{i,j}|$.
Then
\begin{equation*}
   N_{0,0}=N_{1,1}=2^{m-2}+(-1)^{\frac{m^2-1}{8}}2^{\frac{m-3}{2}}, N_{1,0}=N_{0,1}=2^{m-2}-(-1)^{\frac{m^2-1}{8}}2^{\frac{m-3}{2}}.
\end{equation*}

}
\begin{proof}Since $x\mapsto x^{2^s+1}$ is a permutation on $\mathbb{F}_q$, (1) is trivial.

(2) It follows with the fact ${\rm Tr}(x)={\rm Tr}(x^2)$.

(3) By the same reason as (1), we have
\begin{eqnarray*}
  0 &=& C_m^{(s)}(1,0)=\sum_{x\in L}\chi_m(x^{2^s+1})\\
  &=&\sum_{x\in L}\chi_m((x+h)^{2^s+1}) \mbox{ for any $h\in  L$}\\
  &=&\chi_m(h^{2^s+1})\sum_{x\in L}\chi_m(x^{2^s+1}+(h^{2^k}+h^{2^{m-s}})x)\\
  &=&\chi_m(h^{2^s+1})C_m^{(s)}(1,h^{2^s}+h^{2^{m-s}}).
\end{eqnarray*}
Denote
\begin{equation*}
    \mathfrak{H}=\{h^{2^s}+h^{2^{m-s}}|h\in L\}, \mathfrak{T}_0=\{x\in L|{\rm Tr}(x)=0\}.
\end{equation*}
If $h_1^{2^s}+h_1^{2^{m-s}}=h_2^{2^s}+h_2^{2^{m-s}}$, then $h_1^{2^{2k}}+h_1=h_2^{2^{2s}}+h_2$ and thus $h_1-h_2\in \mathbb{F}_{2^{2s}}\cap \mathbb{F}_{2^m}=\mathbb{F}_2$. Therefore, we know that the cardinality of $\mathfrak{H}$ is $2^{m-1}=|\mathfrak{T}_0|$. Moreover, for every $h\in L$, ${\rm Tr}(h^{2^s}+h^{2^{m-s}})=0$. Thus we have
$\mathfrak{H}=\mathfrak{T}_0$ and the conclusion follows.

(4) If ${\rm Tr}(a)=1$, then there exist exactly two elements $h,h+1\in L$ such that $a=h^{2^s}+h^{2^{m-s}}+1$. Thus
\begin{eqnarray*}
  C_m^{(s)}(1,a) &=& \sum_{x\in L}\chi_m(x^{2^s+1}+(h^{2^s}+h^{2^{m-s}}+1)x) \\
&=&\sum_{x\in L}\chi_m(x^{2^s+1}+(h^{2^s}+h^{2^{m-s}}+1)x)\\
&=&\chi_m(h^{2^s+1}+h)\sum_{x\in L}\chi_m(x^{2^s+1}+h^{2^k}x+hx^{2^s}+h^{2^s+1}+x+h)\\
&=&\chi_m(h^{2^s+1}+h)\sum_{x\in L}\chi_m((x+h)^{2^s+1}+(x+h))\\
&=&\chi_m(h^{2^s+1}+h)\sum_{x\in L}\chi_m(x^{2^s+1}+x)\\
&=&\chi_m(h^{2^s+1}+h)C_m^{(s)}(1,1).
\end{eqnarray*}
The desired result now follows with Lemma \ref{LMW}.

(5) Since $\sum_{a\in L}C_m^{(s)}(1,a)=q$, by (3) and (4), we know that
\begin{equation*}
    M_{+}+M_{-}=2^{m-1}, M_{+}-M_{-}=(-1)^{\frac{m^2-1}{8}}2^{\frac{m-1}{2}}.
\end{equation*}
Solving this system of equations leads to the result.

(6) It is obvious that
\begin{eqnarray*}
  N_{i,j}&=&\frac{1}{4}\sum_{h\in \mathbb{F}_{2^m}}\left[1+(-1)^i\chi_m(h^{2^s+1})\right] \left[1+(-1)^j\chi_m(h)\right]\\
  &=&\frac{1}{4}[2^m+(-1)^{i+j}C_m^{(s)}(1,1)].
\end{eqnarray*}
The result is then follows with Lemma \ref{LMW}.
\end{proof}

By Proposition \ref{lem-sb}, one can determine the exponential sum $\sum_{u\in \mathbb{F}_{2^m}^*}\chi_m(\mu(u^{2^s+1}+u))$ and the frequency of it explicitly.

\subsection{The evaluation of $r(l)$}It is evident that
\begin{eqnarray*}
  &&r(l)=2^m+\sum_{a\in \mathbb{F}_{2^n}\setminus \mathbb{F}_{2^m}}\chi_n((a^{2^m}+a)L(a))\\
  &=& 2^m+\sum_{\lambda\in E, u\in \mathbb{F}_{2^m}^*}\chi_n(u L(\lambda u)) \\
  &=&2^m+\sum_{\lambda\in E, u\in \mathbb{F}_{2^m}^*}\chi_m(u L(u(\lambda+\overline{\lambda})))\\
  &=&2^m+\sum_{\lambda\in E, u\in \mathbb{F}_{2^m}^*}\chi_m(u L(u))\\
  &=&2^m\sum_{ u\in \mathbb{F}_{2^m}}\chi_m(u L(u))\\
  &=&2^m\sum_{u\in \mathbb{F}_{2^m}}\chi_m(\sum_{i=0}^k\alpha_iu^{2^{a_i}+1}).
\end{eqnarray*}

We note that for $f(x)=\sum_{i=0}^k\alpha_ix^{2^{a_i}+1}$, the Weil sum $\sum_{u\in \mathbb{F}_{2^m}}\chi_m(f(u))$ has been studied for many authors, see \cite{carlitz, coulter, feng, hou, zhang} etc. In most cases, this Weil sum $\sum_{u\in \mathbb{F}_{2^m}}\chi_m(f(u))$ can be explicitly calculated.
\section*{Acknowledgement}

The authors would like to express their grateful thanks to the referees for their valuable comments and suggestions. The work of this paper is supported by the NUAA Fundamental Research Funds, No. 2013202 and NNSF of China under Grant No. 11371011.

% ----------------------------------------------------------------

\end{document}